\documentclass[11pt,a4paper]{article}
\usepackage[T1]{fontenc}
\usepackage{amsmath,amssymb,amsthm}
\usepackage[margin=29mm]{geometry}
\usepackage{hyperref}
\hypersetup{colorlinks=true,linkcolor=blue,citecolor=blue,urlcolor=blue}
\allowdisplaybreaks
\newtheorem{theorem}{Theorem}
\newtheorem{proposition}{Proposition}[section]
\newtheorem{lemma}[proposition]{Lemma}

\theoremstyle{definition}
\newtheorem{definition}[proposition]{Definition}
\newtheorem{example}[proposition]{Example}
\newtheorem{remark}[proposition]{Remark}
\newcommand{\Z}{\mathbb Z}
\newcommand{\Q}{\mathbb Q}
\newcommand{\R}{\mathbb R}
\newcommand{\Gm}{\mathbb G_m}
\newcommand{\OO}{\mathcal O}
\newcommand{\GU}{\operatorname{GU}}
\newcommand{\GSO}{\operatorname{GSO}}
\newcommand{\GO}{\operatorname{GO}}
\newcommand{\SO}{\operatorname{SO}}
\newcommand{\SU}{\operatorname{SU}}
\newcommand{\UU}{\operatorname U}
\newcommand{\GL}{\operatorname{GL}}
\newcommand{\SL}{\operatorname{SL}}

\newcommand{\Des}{\operatorname{Des}}

\newcommand{\N}{\operatorname N}
\newcommand{\Tr}{\operatorname{Tr}}
\title{Zeta functions of quasi-split reductive algebraic groups over local fields.}
\author{Marcus du Sautoy\thanks{Mathematical Institute, Andrew Wiles Building,
Woodstock Road, Oxford OX2 6GG, UK.}}
\date{}
\begin{document}
\maketitle
\begin{abstract}
We calculate zeta functions of quasi-split reductive algebraic groups over local fields using the classification of Tits and prove
a functional equation. We give explicit examples of these zeta functions for non-split unitary and orthogonal groups.
\end{abstract}
\tableofcontents

\section{Introduction}
Let $K$ be a global field, and let $v$ be a nonarchimedean place.
Write $k=K_v$, $\OO=\OO_v$, $\pi$ for a uniformizer, and
$q=|\OO/\pi\OO|$. Normalize the valuation by $\omega(\pi)=1$ and
the absolute value by $|x|=q^{-\omega(x)}$. Let $G$ be a connected
reductive $K$-group and let $\rho:G\hookrightarrow\GL(V)$ be a
faithful representation, in the sense of a closed immersion. Fix an
integral lattice in $V$ and put $N=\dim_K V$. Define the zeta function of $G$ as
\begin{equation}\label{eq:definition}
 Z_{G,\rho,v}(s)=\int_{G_v^+}|\det\rho(g)|^s\,d\mu_v(g),\qquad
 G_v^+=\{g\in G(k):\rho(g)\in\operatorname{End}_{\OO}(V_{\OO})\},
\end{equation}
where $\mu_v(G(\OO))=1$. 

In this paper we generalise the split-group calculations of Igusa \cite{Igusa}, du Sautoy and
Lubotzky \cite{duSL}, and Berman \cite{Berman} to quasi-split algebraic groups using the classification of Tits \cite{Tits} using the Bruhat decomposition, cone reciprocity and unramified
relative root data. The generalisation to quasi-split algebraic groups is relevant to the calculation of pro-isomorphic zeta functions of nilpotent groups as investigated in \cite{duSL}. In that paper the    algebraic automorphism group of the nilpotent group was assumed to have a split reductive part. This paper allows us to consider non-split examples.

In this paper we establish the finite uniformity of the zeta functions of quasi-split reductive algebraic groups and under additional conditions we prove that the rational functions satisfy a functional equation. The local cone and
its central translation condition are defined in Sections~\ref{sec:local}
and~\ref{sec:reciprocity}.
\begin{theorem}\label{thm:main}
There are finitely many rational functions
$W_1(X,Y),\ldots,W_t(X,Y)\in\Q(X,Y)$ such that, for almost all $v$,
\[
 Z_{G,\rho,v}(s)=W_{i(v)}(q_v,q_v^{-s}).
\]
The integral converges for $\Re(s)$ sufficiently large.
The index $i$ records an unramified local type. For each such type,
choose a representative completion $k_i=K_{v_i}$ and a maximal
$k_i$-split torus $S_i\subseteq G_{k_i}$. Write
$\Phi_i=\Phi(G_{k_i},S_i)$ for its relative root system, choose
positive roots $\Phi_i^+$, and let $m_{i,a}$ be the dimension over
$k_i$ of the relative root space for $a$.
For any type whose chamber cone satisfies the central
translation condition of Definition~\ref{def:translation}, its rational
function satisfies
\begin{equation}\label{eq:mainfe}
 W_i(X^{-1},Y^{-1})=(-1)^{r_i}X^{L_i}Y^{b_i}W_i(X,Y),
\end{equation}
where
\[
 r_i=\dim S_i,\qquad
 L_i=\sum_{a\in\Phi_i^+}m_{i,a},\qquad
 b_i=\langle(\det\rho)|_{S_i},\eta_i\rangle.
\]
Here $\eta_i\in X_*(S_i)$ is an integral cocharacter satisfying
the central translation condition; its meaning and the weighted-length
interpretation of $L_i$ are explained below.
\end{theorem}
The sign involves the local split rank, not the dimension $N$ of the
representation. We use $q_v$, rather than the underlying rational prime;
this also covers general global fields.

The following example illustrates why an extra hypothesis is needed to prove the functional equation.

\begin{example}
\label{ex:counterexample}
Take $G=\Gm^2$ and
\[
 \rho(t,u)=\operatorname{diag}(t,t^2u^3,tu).
\]
This is faithful: its first and third characters generate $X^*(G)$.
For $x=\omega(t)$ and $y=\omega(u)$ the integrality conditions are
$x\geq0$, $2x+3y\geq0$, $x+y\geq0$. The last follows from the
first two. Since $\det\rho=t^4u^4$ and each valuation coset has
measure one, writing $x=3a+r$ gives
\begin{align*}
 Z(s)&=\sum_{x\geq0,\ y\geq\lceil-2x/3\rceil}Y^{4(x+y)}
      =\frac{1+2Y^4}{(1-Y^4)^2},\qquad Y=q^{-s}.
\end{align*}
Indeed the three possible residues $r=0,1,2$ give minimal degrees
$4a,4a+4,4a+4$. Therefore
\[
 \frac{W(Y^{-1})}{W(Y)}
   =\frac{Y^4(Y^4+2)}{1+2Y^4}
\]
is not a signed monomial. Thus an unconditional functional equation for
all connected reductive groups and all faithful representations is false.
Split-group criteria and further counterexamples occur in
\cite{Berman}; the issue is already present before quasi-split descent.
\end{example}

\section{Unramified local data and integral cones}\label{sec:local}
 Assume that $G$ is quasi-split over $k$ and
splits over an unramified extension, and take a hyperspecial subgroup
$\mathcal K=G(\OO)$ compatible with $\rho$. Let $S$ be a maximal
$k$-split torus and $T=Z_G(S)$. Quasi-splitness implies that $T$ is a
maximal torus. Set
\[
 \Lambda=X_*(S),\quad E=\Lambda\otimes_{\Z}\R,\quad
 W=N_G(S)(k)/T(k),\quad r=\dim S.
\]
The valuation quotient $T(k)/T(\OO)$ is identified with $\Lambda$;
for the unramified situation see \cite[\S1.3]{Tits}. We use the
positive-valuation convention: $t_\lambda=\lambda(\pi)$.
Thus our $\lambda$ is the negative of Tits's translation vector
$\nu(t_\lambda)$. This choice makes integrality inequalities nonnegative.

Let $\Phi$ be the relative root system, which may be nonreduced, and
choose positive roots $\Phi^+$ and simple indivisible roots
$\alpha_1,\ldots,\alpha_\ell$. Denote by $m_a$ the dimension over
$k$ of the relative root space for $a$. Define
\begin{equation}\label{eq:AL}
 A=\sum_{a\in\Phi^+}m_a a,\qquad
 L(w)=\sum_{\substack{a\in\Phi^+\\w(a)\in-\Phi^+}}m_a,
 \qquad L=L(w_0)=\sum_{a\in\Phi^+}m_a,
\end{equation}
where $w_0$ is the longest element. Equivalently, $L(w)$ is the
weighted Coxeter length with simple parameter $q^{d_i}$, where
$d_i=L(s_i)$. In a nonreduced simple direction $a$, the parameter
is $q^{m_a+m_{2a}}$. These are the parameters for the finite Weyl group
at the chosen hyperspecial point, rather than all the parameters of
the affine diagram. Put
\[
 D=\{x\in E:\langle\alpha_i,x\rangle\geq0\ (1\leq i\leq\ell)\},
 \qquad \delta=\det\rho|_S.
\]
Let $\Omega_\rho$ be the multiset of weights of $\rho|_S$ and put
\begin{equation}\label{eq:cones}
 C_\rho=\{x\in E:\langle\chi,x\rangle\geq0\ \text{for all }
                 \chi\in\Omega_\rho\},\qquad H=C_\rho\cap D.
\end{equation}
The multiset is $W$-invariant. After choosing the compatible lattice,
\begin{equation}\label{eq:integrality}
 \rho(t_\lambda)\text{ is integral}\ \Longleftrightarrow\
 \lambda\in C_\rho,\qquad
 |\det\rho(t_\lambda)|^s=q^{-s\langle\delta,\lambda\rangle}.
\end{equation}

\begin{lemma}\label{lem:weights}
The weights of a faithful representation of a split torus $S$ generate
$X^*(S)$ as an abelian group. Consequently $C_\rho$ is pointed and
$\langle\delta,x\rangle>0$ for every nonzero $x\in C_\rho$.
\end{lemma}
\begin{proof}
The scheme-theoretic kernel of the weight homomorphism
$S\to\Gm^{\Omega_\rho}$ is the diagonalizable group whose character
group is $X^*(S)/\langle\Omega_\rho\rangle$. Faithfulness makes this
quotient zero. If both $x$ and $-x$ are in $C_\rho$, every weight
vanishes on $x$, so $x=0$. On $C_\rho$, the value of $\delta$ is the
sum of the nonnegative weight values, with their multiplicities; it
can be zero only at zero.
\end{proof}
The scheme-theoretic qualification is needed in positive characteristic;
triviality of the kernel on rational points alone is insufficient.

\subsection{Affine roots and their multiplicities}
We recall the portion of the Bruhat--Tits notation needed for explicit
calculations. An affine root is represented by $a+c$, with $a$ a
relative root and $c$ a jump in its valued root-group filtration. If
$U_{a,c+}=\bigcup_{t>c}U_{a,t}$, the number $d(a+c)$ is the
residue-field dimension of the associated graded quotient after removing
the contribution of $U_{2a}$. For an unramified group, at a compatible
hyperspecial origin the jumps can be taken to be integral and
\begin{equation}\label{eq:periodic}
 d(a+j)=m_a\qquad(j\in\Z).
\end{equation}
This follows by splitting over an unramified extension: the relative
root space is the descent of the sum of absolute root spaces restricting
to $a$, and its integral lattice has graded pieces of constant dimension.
For multipliable roots this applies separately to the quotient by
$U_{2a}$ and to $U_{2a}$. This is precisely the distinction in
\cite[\S\S1.4--1.8, 1.15]{Tits} and \cite[\S4.2]{BT2}.
It is not an assertion for arbitrary ramified root data.

If $a$ and $2a$ are roots, their affine functions can have the same
wall; a wall at $a+j=0$ then has total exponent $m_a+m_{2a}$,
whereas a wall at $2a+(2j+1)=0$ can have exponent $m_{2a}$ alone.
Both must be counted. Also, for a negative pairing, the count of jumps
in $0<c\leq\langle a,\lambda\rangle$ is zero, not that negative
pairing. Linear jump-count formulas are used only after moving into
the corresponding positive chamber.

There is a small notational point concerning an \emph{\'echelonnage}.
If affine roots are regarded as half-spaces, their defining functions
are determined only up to positive scalar. Thus the relative root
paired with a half-space is a positive multiple of the linear normal
of a defining affine function. If the representative is specifically
$a+c$, its linear part is exactly $a$. The two projections of the
\'echelonnage onto the relative roots and the affine half-spaces are
surjective, with the factors ordered accordingly. For the present
unramified calculation, the filtration dimensions above supply the
required information without a classification of all ramified
\'echelonnages.

\section{Double cosets and the generating-function formula}
For $J\subseteq\{1,\ldots,\ell\}$ let $W_J=\langle s_j:j\in J\rangle$,
and put
\[
 P_J(t)=\sum_{w\in W_J}t^{L(w)},\qquad P_W(t)=\sum_{w\in W}t^{L(w)}.
\]
Write $\Des(w)=\{i:w(\alpha_i)\in-\Phi^+\}$ for right descents,
and let $J(\lambda)=\{i:\langle\alpha_i,\lambda\rangle=0\}$.

\begin{proposition}\label{prop:cartan}
The hyperspecial Cartan decomposition and its volume formula are
\begin{align}
 G(k)&=\coprod_{\lambda\in\Lambda\cap D}\mathcal Kt_\lambda\mathcal K,
                                                       \label{eq:cartan}\\
 \mu(\mathcal Kt_\lambda\mathcal K)
  &=q^{\langle A,\lambda\rangle}
        \frac{P_W(q^{-1})}{P_{J(\lambda)}(q^{-1})}.        \label{eq:volume}
\end{align}
Consequently, if
\[
 H_w=\{\lambda\in\Lambda\cap H:
       \langle\alpha_i,\lambda\rangle>0\text{ for every }i\in\Des(w)\},
\]
then
\begin{equation}\label{eq:weylsum}
 Z_{G,\rho,v}(s)=
 \sum_{w\in W}q^{-L(w)}
       \sum_{\lambda\in H_w}
          q^{\langle A,\lambda\rangle-s\langle\delta,\lambda\rangle}.
\end{equation}
\end{proposition}
\begin{proof}
The Cartan and Iwahori decompositions follow from the valued root datum;
see \cite[\S3.3]{Tits}, \cite[\S4.4]{BT1}, and the split predecessor
\cite{IM}. Choose an Iwahori
$\mathcal I\subset\mathcal K$. Its finite Bruhat cells have volumes
$q^{L(w)}\mu(\mathcal I)$, so
$\mu(\mathcal I)P_W(q)=1$.
For a more explicit volume computation, put
$\mathcal K_\lambda=\mathcal K\cap
t_\lambda\mathcal Kt_\lambda^{-1}$ and $J=J(\lambda)$.
Its reduction is a parabolic subgroup with Weyl group $W_J$ in the
reductive special fibre. Write $L_J$ for the weighted longest length
in $W_J$. The finite flag variety has
\[
 \sum_{w\in W^J}q^{L(w)}
 =q^{L-L_J}\frac{P_W(q^{-1})}{P_J(q^{-1})}
\]
points, where $W^J$ denotes minimal right coset representatives.
The identity follows from length-additive parabolic factorization and
the longest-element symmetry in $W$ and $W_J$.
The root directions with $\langle a,\lambda\rangle>0$ contribute
further pro-unipotent index
\[
 q^{\sum_{a\in\Phi^+,\ \langle a,\lambda\rangle>0}
                 m_a(\langle a,\lambda\rangle-1)}
       =q^{\langle A,\lambda\rangle-(L-L_J)}.
\]
Indeed the first filtration step is accounted for by reduction to
the parabolic, and each further step contributes $q^{m_a}$;
ordered root filtrations give the product also in nonreduced directions.
Multiplying these two indices proves \eqref{eq:volume}, since
$\mu(\mathcal Kt_\lambda\mathcal K)=[\mathcal K:\mathcal K_\lambda]$.
Equivalently the finite correction to $q^{\langle A,\lambda\rangle}$ is
\[
 \sum_{w\in W^{J(\lambda)}}q^{-L(w)}
       =P_W(q^{-1})/P_{J(\lambda)}(q^{-1}).
\]
At $\lambda=0$ this correction is one, as required by normalization.
Both integrality and the determinant norm are constant on each
$\mathcal K$-double coset. Finally,
$w\in W^J$ is equivalent to $\Des(w)\cap J=\varnothing$.
Expanding the quotient and interchanging the two sums gives
\eqref{eq:weylsum}.
\end{proof}
In particular there is no separate factor $\mu(\mathcal I)$ in
\eqref{eq:weylsum}. The compact normalization has already cancelled it.
For groups with split center the group of length-zero extended affine
Weyl elements can be infinite. The decomposition above does not require
it to have finite index or a point stabilizer in the reduced building
to be compact; we use the hyperspecial integral subgroup throughout.

\subsection{Rationality and the passage to global fields}
For a fixed local type define
\begin{equation}\label{eq:Wi}
 W(X,Y)=\sum_{w\in W}X^{-L(w)}
     \sum_{\lambda\in H_w}
        X^{\langle A,\lambda\rangle}Y^{\langle\delta,\lambda\rangle}.
\end{equation}
Each inner series is rational. Indeed a rational pointed cone has a
finite triangulation by rational simplicial cones. For a simplicial cone
with independent integral ray vectors $v_1,\ldots,v_d$, the generating
function is a finite numerator from a fundamental parallelepiped divided
by
\[
 \prod_{j=1}^d
    (1-X^{\langle A,v_j\rangle}Y^{\langle\delta,v_j\rangle}).
\]
Deleting selected rational faces or using half-open pieces preserves
rationality, by inclusion-exclusion. This argument also applies when
$H$ has dimension smaller than $r$, using its induced lattice.
By Lemma~\ref{lem:weights}, every nonzero ray has positive determinant
degree. Thus the specialization converges when $\Re(s)$ exceeds the
finitely many ratios $\langle A,v_j\rangle/\langle\delta,v_j\rangle$.

Choose a finite Galois extension of $K$ over which $G$ splits. Outside
finitely many places it is unramified, $G$ has a reductive smooth
integral model, and the representation and lattice are compatible with
that model. The special fibre is a connected reductive group over a
finite field, hence has a rational Borel by Lang's theorem. The smooth
scheme of Borel subgroups lifts this point over the henselian valuation
ring. Consequently the local group is quasi-split and its integral
subgroup is hyperspecial; see also \cite[\S3.9]{Tits}.
This uses a finite splitting extension, not an assertion that the whole
separable closure is unramified at almost all places.

There are only finitely many actions of Frobenius arising from this fixed
finite Galois group on the absolute root datum and on the representation
weights. They determine the lattices, relative multiplicities, cones,
and lengths in \eqref{eq:Wi}, up to isomorphism. Choosing one rational
function for each occurring type proves the rationality assertion of
Theorem~\ref{thm:main}.

\section{Reciprocity and the functional equation}\label{sec:reciprocity}
Assume in this section that $H$ has full dimension $r$ and that its
distinct facets consist of the simple root walls and some further walls.
Write an irredundant description as
\begin{equation}\label{eq:facets}
 H=\{x:\alpha_i(x)\geq0\ (1\leq i\leq\ell),\quad
                  f_j(x)\geq0\ (1\leq j\leq h)\},
\end{equation}
where each $f_j$ is integral and primitive on $\Lambda$.
No primitive normalization of the $\alpha_i$ is needed since only
their signs are used.
\begin{definition}[Central translation condition]\label{def:translation}
The chamber cone satisfies the central translation condition if it has
the full-dimensional facet description \eqref{eq:facets} and there is
$\eta\in\Lambda$ such that
\begin{equation}\label{eq:eta}
 \alpha_i(\eta)=0\quad(1\leq i\leq\ell),\qquad
 f_j(\eta)=1\quad(1\leq j\leq h).
\end{equation}
\end{definition}
This is an explicit sufficient condition, not a claim of necessity for
every univariate specialization. The equations can be checked by
integer linear algebra once the primitive facet forms are known.

For $I\subseteq\{1,\ldots,\ell\}$ put
\[
 H_I^{\Z}=\{\lambda\in\Lambda\cap H:
                         \alpha_i(\lambda)>0\ (i\in I)\},\quad
 F_I(\mathbf z)=\sum_{\lambda\in H_I^{\Z}}\mathbf z^\lambda.
\]
The symbols denote rational functions when variables are inverted.
\begin{lemma}\label{lem:reciprocity}
Under the central translation condition,
\begin{equation}\label{eq:reciprocity}
 F_I(\mathbf z^{-1})=(-1)^r\mathbf z^\eta F_{I^c}(\mathbf z).
\end{equation}
\end{lemma}
\begin{proof}
Choose $u\in E$ with $\alpha_i(u)<0$ for $i\in I$ and
$\alpha_i(u)>0$ otherwise. Adding a sufficiently large positive
multiple of $\eta$ makes every $f_j(u)>0$ without changing these
signs. The lattice points $\lambda$ with
$\lambda+\epsilon u\in H$, for sufficiently small $\epsilon>0$,
are precisely $H_I^{\Z}$. Uniform smallness is possible because the
facet forms are rational on $\Lambda$.

The half-open cone reciprocity identity gives
\[
 F_I(\mathbf z^{-1})=(-1)^r
 \sum_{\substack{\lambda\in\Lambda:\ f_j(\lambda)>0\ \forall j,\\
                \alpha_i(\lambda)>0\ (i\notin I),\
                \alpha_i(\lambda)\geq0\ (i\in I)}}\mathbf z^\lambda.
\]
For completeness this version follows by triangulating $H$ into
simplicial cones and using a generic perturbation of $u$ to assign
their shared faces. On each piece, inversion of the ray denominators
contributes $(-1)^r$ and replaces the chosen half-open fundamental
parallelepiped by the complementary one. The internal faces cancel
on summing, leaving the displayed complementary facet conditions.
This is the usual cone reciprocity argument; see \cite{Stanley,BeckSottile}.
Primitivity makes $f_j(\lambda)>0$ equivalent to
$f_j(\lambda)\geq1$. Translation by $\eta$ therefore identifies
the last indexing set with $\eta+H_{I^c}^{\Z}$, proving the claim.
\end{proof}

\begin{proof}[Completion of the proof of Theorem~\ref{thm:main}]
It remains to prove the functional equation. Specialize
$\mathbf z^\lambda$ to
$X^{\langle A,\lambda\rangle}Y^{\langle\delta,\lambda\rangle}$
in Lemma~\ref{lem:reciprocity}. Since $A(\eta)=0$, the translation
factor is $Y^b$, where $b=\delta(\eta)$. The two elementary
weighted Weyl-group identities are
\[
 \Des(w_0w)=\Des(w)^c,\qquad L(w_0w)=L-L(w).
\]
They follow respectively by reversal of positive roots under $w_0$
and by the complementary inversion sets, with their root multiplicities.
Consequently \eqref{eq:Wi} gives
\begin{align*}
 W(X^{-1},Y^{-1})
  &=(-1)^rY^b\sum_{w\in W}X^{L(w)}F_{\Des(w)^c}(X,Y)\\
  &=(-1)^rX^LY^b\sum_{u\in W}X^{-L(u)}F_{\Des(u)}(X,Y)\\
  &=(-1)^rX^LY^bW(X,Y).
\end{align*}
Together with the preceding rationality and global arguments this proves
all assertions. There is no application of reciprocity to an arbitrary
cone without checking the translation condition.
\end{proof}

\section{Unitary example and the calculation for \texorpdfstring{$\GU(3)$}{GU(3)}}
\label{sec:GU3}
Let $E/k$ be a quadratic extension, with conjugation $x\mapsto\bar x$.
For odd dimension $n=2r+1$, take the quasi-split hermitian form
\begin{equation}\label{eq:hermitian}
 h(x)=\sum_{i=1}^r(\bar x_{-i}x_i+\bar x_ix_{-i})+\bar x_0x_0.
\end{equation}
For even dimension use the same hyperbolic sum, omitting $x_0$.
We make the group and representation conventions explicit:
\[
 \UU(h)=\{g:h(gx)=h(x)\},\qquad
 \SU(h)=\UU(h)\cap\SL_n(E),
\]
\begin{equation}\label{eq:GU}
 \GU(h)=\{g:h(gx)=c(g)h(x),\ c(g)\in k^\times\}.
\end{equation}
Here $\GU$ denotes similitudes. In \cite{Tits} Tits's Example~1.15 concerns $\SU$.
It supplies the derived-group root data needed for \eqref{eq:GU}.
Our faithful $k$-representation $\rho$ is the action on the underlying
$2n$-dimensional $k$-space of $E^n$; the integral lattice is $\OO_E^n$.
Taking determinants of $\bar g^{\,t}Jg=c(g)J$ yields
\begin{equation}\label{eq:detnorm}
 \det_k\rho(g)=\N_{E/k}(\det_Eg)=c(g)^n.
\end{equation}
Thus the integrand is $|c(g)|^{ns}$. This specifies the normalization
of $s$; using a square root of this norm would rescale $s$.

\subsection{The unramified nonsplit root groups}
Assume now that $E/k$ is unramified. In dimension three, a positive
root group in the ordered coordinates $(-1,0,1)$ is
\begin{equation}\label{eq:rootgroup}
 u(c,d)=\begin{pmatrix}1&-\bar c&d\\0&1&c\\0&0&1\end{pmatrix},
 \qquad \bar c c+d+\bar d=0.
\end{equation}
Its center is $u(0,d)$ with $d+\bar d=0$. Tits's filtration is
$\phi_a(u(c,d))=\frac12\omega(d)$ and
$\phi_{2a}(u(0,d))=\omega(d)$. The trace map of the unramified
extension is surjective on integers, and the trace-zero lattice has
rank one over $\OO$. Hence, as recorded in
\cite[Example~1.15, pp.~40--42]{Tits},
\begin{equation}\label{eq:SU3d}
 \Gamma'_a=\Gamma'_{2a}=\Z,\qquad
 d(a+j)=2,\quad d(2a+j)=1\quad(j\in\Z).
\end{equation}
For example, modulo $\pi$ the group in \eqref{eq:rootgroup} has
$q^3$ elements, its center has $q$ elements, and its quotient has
$q^2$ elements. In particular the two dimensions in
\eqref{eq:SU3d} must not be interchanged.

The affine walls alternate between total dimensions $3=2+1$ and $1$.
One may depict the rank-one affine diagram, with its fat bond, as
\[
 \underset{3}{\overset{\times}{\circ}}
 \mathrel{\blacktriangleleft\!\rule[0.45ex]{1.4cm}{1.2pt}}
 \underset{1}{\circ}.
\]
Its building is a tree with vertex degrees $q^3+1$ and $q+1$.
Choose the hyperspecial vertex of degree $q^3+1$.
The finite reflection there has parameter $q^3$, whereas a translation
by one cocharacter crosses walls of total dimension $4$.

\subsection{The unramified nonsplit factor}
Every dominant Cartan representative of $\GU(3)(k)$ can be chosen as
\[
 t_{m,a}=\operatorname{diag}(\pi^{m+a},\pi^m,\pi^{m-a}),
       \qquad m\in\Z,\ a\in\Z_{\geq0}.
\]
Its multiplier is $\pi^{2m}$. That the multiplier valuation is even
also follows from \eqref{eq:detnorm}: the valuation of a norm from an
unramified quadratic extension is even, and $3\omega(c)$ is even.
The parameters $(m,a)$ give the full valuation lattice; compact norm-one
torus factors do not contribute additional indices.
Integrality is equivalent to $m\geq a$, and the determinant valuation
over $k$ is $6m$. By \eqref{eq:SU3d} and
Proposition~\ref{prop:cartan},
\begin{equation}\label{eq:SU3vol}
 \mu(\mathcal Kt_{m,a}\mathcal K)=
 \begin{cases}
 1,&a=0,\\
 (1+q^{-3})q^{4a},&a\geq1.
 \end{cases}
\end{equation}
The tree gives an independent interpretation: at distance $2a$ from a
vertex of degree $q^3+1$, the number of vertices of the same type is
$(q^3+1)q^{4a-3}$.

Put $T=q^{-6s}$. Summing first over $m\geq a$ gives
\begin{align}
 Z_{\GU(3),\rho}(s)
  &=\frac1{1-T}\left(1+(1+q^{-3})
                              \sum_{a\geq1}(q^4T)^a\right)\notag\\
  &=\boxed{\frac{1+q^{1-6s}}
               {(1-q^{-6s})(1-q^{4-6s})}}.             \label{eq:GU3inert}
\end{align}
Its abscissa of convergence is $2/3$. As a rational function,
\begin{equation}\label{eq:GU3FE}
 W_{\rm ns}(X^{-1},Y^{-1})=X^3Y^6W_{\rm ns}(X,Y).
\end{equation}
Here the split rank is two, $L=3$, and the central translation
$(m,a)=(1,0)$ has determinant degree six. The first coefficients in
$T$ are $1$, $q^4+q+1$, and $q^8+q^5+q^4+q+1$.

\subsection{The split factor}
If the quadratic algebra splits, $E=k\times k$, then
\[
 \GU(n)\simeq\GL_n\times\Gm,\qquad
 \rho(g,c)=(g,c(g^{-1})^t).
\]
For $n=3$ a dominant valuation tuple is
$a_1\geq a_2\geq a_3$, together with $b=\omega(c)$.
Integrality requires $0\leq a_3\leq a_2\leq a_1\leq b$.
Set
\[
 u_0=b-a_1,\quad u_1=a_1-a_2,\quad
 u_2=a_2-a_3,\quad u_3=a_3.
\]
These are independent nonnegative integers. The determinant degree is
$3\sum u_i$, and $A=2u_1+2u_2$. The six permutations of $S_3$
in \eqref{eq:weylsum} give, with $T=q^{-3s}$,
\begin{equation}\label{eq:GU3split}
 \boxed{Z_{\GU(3),\rho}(s)=
 \frac{1+2(q+1)T+qT^2}{(1-T)^2(1-q^2T)^2}}.
\end{equation}
It too has abscissa $2/3$ and satisfies
$W_{\rm sp}(X^{-1},Y^{-1})=X^3Y^6W_{\rm sp}(X,Y)$.
At a good global place the factor is \eqref{eq:GU3inert} or
\eqref{eq:GU3split} according to whether the quadratic extension is
inert or split. Ramified places are outside these formulas.

\begin{remark}[Isometries rather than similitudes]
For the same natural representation of $\UU(n)$ or $\SU(n)$, the
determinant over $k$ has absolute value one. An integral group matrix
then has integral inverse, so $G^+=G(\OO)$ and the normalized zeta
function is $1$. Thus the nontrivial formulas above concern
similitudes and do not describe the natural representation of the
special unitary group in Tits's example itself.
\end{remark}

\section{Finite formulas for \texorpdfstring{$\GU(n)$}{GU(n)} in every dimension}\label{sec:GUn}
We keep the representation and determinant convention of
\eqref{eq:detnorm}. All the formulas below are finite rational
expressions; their numerators can be calculated by enumerating a finite
Weyl group. They apply to the quasi-split form and its self-dual
hyperspecial lattice at an unramified nonsplit place, or to the split
quadratic algebra at a split place.

\subsection{Odd dimension at an unramified nonsplit place}
Let $n=2r+1$. The relative roots of the derived group are
\[
 \{\pm e_i\pm e_j:i<j\}\ \cup\ \{\pm e_i,\pm2e_i:1\leq i\leq r\},
\]
with multiplicities $2$, $2$, and $1$, respectively. The Weyl group
$B_r$ is the signed permutation group, with simple reflections
$s_i$ interchanging $e_i,e_{i+1}$ for $i<r$ and $s_r$ changing the
sign of $e_r$. Its weighted length has parameters
\[
 d_i=2\ (i<r),\quad d_r=3,\quad L(w_0)=r(2r+1).
\]
Define
\begin{equation}\label{eq:oddMi}
 M_i=X^{2i(n-i)}Y^{2n}\qquad(1\leq i\leq r).
\end{equation}
\begin{theorem}\label{thm:odd}
For $n=2r+1$ at an unramified nonsplit place,
\begin{equation}\label{eq:oddformula}
 W_{n,\rm ns}(X,Y)=
 \frac{\displaystyle\sum_{w\in B_r}X^{-L(w)}
                         \prod_{i\in\Des(w)}M_i}
 {(1-Y^{2n})\displaystyle\prod_{i=1}^r(1-M_i)}.
\end{equation}
Moreover
\begin{equation}\label{eq:oddfe}
 W_{n,\rm ns}(X^{-1},Y^{-1})
 =(-1)^{r+1}X^{r(2r+1)}Y^{2n}W_{n,\rm ns}(X,Y).
\end{equation}
\end{theorem}
\begin{proof}
Write the dominant diagonal valuations as
\[
 m+a_1,\ldots,m+a_r,m,m-a_r,\ldots,m-a_1,
 \qquad a_1\geq\cdots\geq a_r\geq0.
\]
Here $m,a_i\in\Z$, integrality is $m\geq a_1$, and the determinant
degree is $2nm$. Tits's multiplicities give
\[
 A(\lambda)=\sum_{i=1}^r4(r-i+1)a_i.
\]
Put $a_{r+1}=0$, $u_0=m-a_1$, and $u_i=a_i-a_{i+1}$.
Then all $u_i$ are independent nonnegative integers and
\[
 m=\sum_{i=0}^ru_i,\qquad
 A(\lambda)=\sum_{i=1}^r2i(n-i)u_i.
\]
A descent in direction $i$ makes $u_i$ positive. Each inner sum
in \eqref{eq:weylsum} is therefore the product of the geometric
series in \eqref{eq:oddformula}. The only non-root facet of $H$
is $u_0=0$; its primitive form takes value one at the central
cocharacter $(m;a_1,\ldots,a_r)=(1;0,\ldots,0)$.
Theorem~\ref{thm:main} proves \eqref{eq:oddfe}.
\end{proof}
The empty-group convention for $r=0$ gives $W_{1,\rm ns}=1/(1-Y^2)$.
The case $r=1$ is \eqref{eq:GU3inert}.

\subsection{Even dimension at an unramified nonsplit place}
Let $n=2r$. The relative roots are $\pm e_i\pm e_j$ of multiplicity
two and $\pm2e_i$ of multiplicity one. Use the same signed
permutation group, now with parameters
$d_i=2$ for $i<r$ and $d_r=1$. Thus $L(w_0)=r(2r-1)$.
Set
\begin{equation}\label{eq:evenMi}
 M_i=X^{2i(n-i)}Y^{2n}\ (i<r),\qquad
 M_r=X^{r^2}Y^n.
\end{equation}
\begin{theorem}\label{thm:even}
For $n=2r$ at an unramified nonsplit place,
\begin{equation}\label{eq:evenformula}
 W_{n,\rm ns}(X,Y)=
 \frac{\displaystyle\sum_{w\in B_r}X^{-L(w)}
                         \prod_{i\in\Des(w)}M_i}
 {(1-Y^{2n})\displaystyle\prod_{i=1}^r(1-M_i)},
\end{equation}
and
\begin{equation}\label{eq:evenfe}
 W_{n,\rm ns}(X^{-1},Y^{-1})
 =(-1)^{r+1}X^{r(2r-1)}Y^{2n}W_{n,\rm ns}(X,Y).
\end{equation}
\end{theorem}
\begin{proof}
Let $b$ be the multiplier valuation. The diagonal valuations can be
taken to be $x_1,\ldots,x_r,b-x_r,\ldots,b-x_1$, where
$x_1\geq\cdots\geq x_r\geq b/2$. All $b,x_i$ are integers;
unlike odd dimension, $b$ need not be even. Introduce
\[
 u_0=b-x_1,\qquad u_i=x_i-x_{i+1}\ (1\leq i<r),\qquad
 u_r=2x_r-b.
\]
These are independent nonnegative integral coordinates: their inverse
is determined by
\[
 b=2\sum_{i=0}^{r-1}u_i+u_r,\qquad
 x_r=\sum_{i=0}^{r-1}u_i+u_r,\qquad
 x_i=x_r+\sum_{j=i}^{r-1}u_j.
\]
In particular no even-multiplier condition is to be imposed.
The determinant degree is $nb$ and, writing $a_i=x_i-b/2$,
\[
 A(\lambda)=\sum_{i=1}^r(4(r-i)+2)a_i
 =\sum_{i=1}^{r-1}2i(n-i)u_i+r^2u_r.
\]
Descent in direction $i$ requires $u_i>0$, yielding
\eqref{eq:evenformula}. The central cocharacter with $b=2$ and
all $x_i=1$ has $u_0=1$ and all other $u_i=0$, proving the
translation condition and \eqref{eq:evenfe}.
\end{proof}
For example, with $T=Y^2$, the case $n=2$ simplifies to
\[
 W_{2,\rm ns}(X,Y)=\frac{1+T}{(1-T^2)(1-XT)}
                 =\frac1{(1-T)(1-XT)}.
\]

\subsection{All dimensions at split places}
Let $\ell(w)$ be ordinary inversion length in $S_n$, and set
$M_i=X^{i(n-i)}Y^n$ for $1\leq i<n$.
\begin{theorem}\label{thm:split}
At a split place,
\begin{equation}\label{eq:splitformula}
 W_{n,\rm sp}(X,Y)=
 \frac{\displaystyle\sum_{w\in S_n}X^{-\ell(w)}
                            \prod_{i\in\Des(w)}M_i}
 {(1-Y^n)^2\displaystyle\prod_{i=1}^{n-1}(1-M_i)},
\end{equation}
and
\begin{equation}\label{eq:splitfe}
 W_{n,\rm sp}(X^{-1},Y^{-1})
 =(-1)^{n+1}X^{n(n-1)/2}Y^{2n}W_{n,\rm sp}(X,Y).
\end{equation}
\end{theorem}
\begin{proof}
As in the split calculation for dimension three, put
$u_0=b-a_1$, $u_i=a_i-a_{i+1}$ for $1\leq i<n$, and
$u_n=a_n$. They independently range over $\Z_{\geq0}$ and
\[
 b=\sum_{i=0}^nu_i,\qquad
 A(\lambda)=\sum_{i=1}^{n-1}i(n-i)u_i,\qquad
 \delta(\lambda)=nb.
\]
The two non-root facets are $u_0=0$ and $u_n=0$. The central
cocharacter $a_i=1$, $b=2$ takes value one on both primitive
forms. Geometric summation proves \eqref{eq:splitformula} and
Theorem~\ref{thm:main} proves \eqref{eq:splitfe}.
\end{proof}
In particular the local sign and the power of $X$ can depend on the
splitting type; no common exponent independent of local type is
asserted in general. The formulas use finite groups of orders
$2^rr!$ and $n!$, so they are directly implementable without an
infinite lattice enumeration.

\subsection{Two further explicit nonsplit examples}
In dimension four, put $T=Y^4$. Enumerating the eight signed
permutations gives numerator
$(1+T)(1+XT)(1+X^3T)$ in \eqref{eq:evenformula}.
Cancelling $1+T$ against $1-T^2$ yields
\begin{equation}\label{eq:GU4}
 W_{4,\rm ns}(X,Y)=
 \frac{(1+XT)(1+X^3T)}{(1-T)(1-X^4T)(1-X^6T^2)},\qquad T=Y^4.
\end{equation}
In dimension five, with $T=Y^{10}$, the eight elements in
\eqref{eq:oddformula} give
\begin{equation}\label{eq:GU5}
 W_{5,\rm ns}(X,Y)=
 \frac{1+(X+X^3+X^4+X^6+X^7+X^9)T+X^{10}T^2}
      {(1-T)(1-X^8T)(1-X^{12}T)}.
\end{equation}
Their functional-equation factors are $-X^6Y^8$ and $-X^{10}Y^{10}$,
respectively.

\section{Nonsplit orthogonal groups}
\label{sec:orthogonal}
Let $E/k$ be an unramified quadratic field extension, and in this
section assume that the residue characteristic is odd. On the
$d=2r+2$ dimensional $k$-space $k^r\oplus E\oplus k^r$, take
\begin{equation}\label{eq:orthogonalform}
 Q_r(x)=\sum_{i=1}^r x_{-i}x_i+\N_{E/k}(x_0),
 \qquad x_0\in E,\quad x_{\pm i}\in k.
\end{equation}
This is precisely the form in \cite[Example~1.16, pp.~42--43]{Tits},
with $r$ in place of Tits's $n$. Its Witt index is $r$ and its
anisotropic kernel is the norm plane, so the even orthogonal group
is quasi-split but nonsplit. Tits assumes $r\geq2$; our main
examples therefore have $d=6,8,10$. We use the natural
$d$-dimensional representation on the lattice
$\OO^r\oplus\OO_E\oplus\OO^r$.

\subsection{Isometries and connected similitudes}
Tits writes $O(Q_r)$ for the full, generally disconnected isometry
group. In its natural representation the integral is
\begin{equation}\label{eq:orthogonalisometry}
 Z_{O(Q_r),\mathrm{nat}}(s)=Z_{\SO(Q_r),\mathrm{nat}}(s)=1.
\end{equation}
Indeed an isometry has unit determinant. An integral isometry thus
has integral inverse and belongs to the integral compact subgroup;
the integrand is one, and that subgroup has measure one.

As in the passage from unitary isometries to $\GU$, the nontrivial
analogue is the connected similitude group
\begin{equation}\label{eq:GSOdefinition}
 G_d^-:=\GSO(Q_r)=
 \{g:Q_r(gx)=c(g)Q_r(x),\quad\det g=c(g)^{r+1}\}.
\end{equation}
Here the superscript minus indicates the nonsplit even orthogonal
form, not a choice of sign for the similitude multiplier. The group
in \eqref{eq:GSOdefinition} is connected reductive, so
Theorem~\ref{thm:main} applies. We normalize Haar measure on its
natural hyperspecial subgroup and write
\[
 Z_{G_d^-,\mathrm{nat}}(s)=W_d^-(q,q^{-s}).
\]

The multiplier of any similitude lies in $\N_{E/k}(E^\times)$.
To see this, $Q_r$ and $cQ_r$ are isometric; their hyperbolic
summands are isometric, and Witt cancellation shows that the norm
plane and its scalar multiple by $c$ are isometric. Evaluating that
plane isometry on a vector of norm one shows that $c$ is a norm.
Conversely, for $z\in E^\times$, multiplication by $z$ on the norm
plane, together with paired scalings on the hyperbolic summands,
realizes the multiplier $\N(z)$ in $\GSO(Q_r)$.
Since the extension is unramified, multiplier valuations are therefore
even. This parity condition is part of the valuation lattice and
must not be discarded.

If one uses the full similitude group $\GO(Q_r)$ instead, its
normalized factor is the same as for $\GSO(Q_r)$. Conjugation on
the norm plane, acting trivially on the hyperbolic summands, is an
integral isometry of determinant minus one. Multiplication by this
element identifies the two components, preserves the integral domain
and determinant norm, and makes their equal contributions add to the
connected-group factor after compact normalization.

\subsection{Relative roots, the lattice, and the general formula}
The relative root system of $\SO(Q_r)$ is reduced of type $B_r$:
\begin{equation}\label{eq:orthogonalroots}
 \Phi=\{\pm e_i\pm e_j:i<j\}\ \cup\ \{\pm e_i:1\leq i\leq r\},
 \qquad m_{e_i\pm e_j}=1,\quad m_{e_i}=2.
\end{equation}
These are the multiplicities encoded by the unramified diagram in
Tits's Example~1.16. One can also obtain them directly from the
absolute root system $D_{r+1}$: the last coordinate restricts to zero
on the split torus, so each pair $e_i+e_{r+1},e_i-e_{r+1}$ restricts
to the same short relative root. The long relative roots remain
single absolute roots. At a hyperspecial origin all jumps are
integral. Thus the signed permutation Weyl group $B_r$ has simple
parameters
\begin{equation}\label{eq:orthogonalparameters}
 d_i=1\ (i<r),\qquad d_r=2,\qquad L(w_0)=r(r+1).
\end{equation}
We keep the conventions $\alpha_i=e_i-e_{i+1}$ for $i<r$ and
$\alpha_r=e_r$, and use right descents as before. In particular
the parameters differ from both unitary cases.

The full dominant valuation lattice is parametrized by
\[
 (m;a_1,\ldots,a_r),\qquad
 m\in\Z,\quad a_1\geq\cdots\geq a_r\geq0,\quad a_i\in\Z.
\]
A representative scales the paired hyperbolic coordinates by
$\pi^{m+a_i}$ and $\pi^{m-a_i}$ and acts on $E$ by multiplication
by $\pi^m$. Its multiplier is $\pi^{2m}$ and its determinant is
$\pi^{dm}$. For completeness, a centralizer element acts by
$z\in E^\times$ on the norm plane and by $t_i,\N(z)t_i^{-1}$ on
the paired lines. Setting $m=\omega_E(z)$ and
$a_i=\omega(t_i)-m$ proves that these integral parameters exhaust
the lattice. The Weyl group orders the $a_i$ and changes their signs;
compact torus factors are absorbed in the hyperspecial subgroup.

Integrality is exactly $m\geq a_1$. The translation exponent and
the determinant degree are
\begin{equation}\label{eq:orthogonalexponents}
 A(\lambda)=\sum_{i=1}^r2(r-i+1)a_i,\qquad
 \delta(\lambda)=dm.
\end{equation}
Put $a_{r+1}=0$ and introduce the independent nonnegative coordinates
\[
 u_0=m-a_1,\qquad u_i=a_i-a_{i+1}\ (1\leq i\leq r).
\]
Then $m=\sum_{i=0}^ru_i$ and
$A(\lambda)=\sum_{i=1}^r i(2r+1-i)u_i$.

\begin{theorem}\label{thm:orthogonal}
Set $d=2r+2$ and
$M_i=X^{i(2r+1-i)}Y^d$ for $1\leq i\leq r$.
For the natural representation of the unramified nonsplit group
$G_d^-=\GSO(Q_r)$,
\begin{equation}\label{eq:orthogonalformula}
 W_d^-(X,Y)=
 \frac{\displaystyle\sum_{w\in B_r}X^{-L(w)}
                          \prod_{i\in\Des(w)}M_i}
 {(1-Y^d)\displaystyle\prod_{i=1}^r(1-M_i)},
\end{equation}
where $L$ uses the parameters \eqref{eq:orthogonalparameters}.
The functional equation is
\begin{equation}\label{eq:orthogonalfe}
 W_d^-(X^{-1},Y^{-1})
  =(-1)^{r+1}X^{r(r+1)}Y^dW_d^-(X,Y).
\end{equation}
\end{theorem}
\begin{proof}
In Proposition~\ref{prop:cartan}, descent in direction $i$ requires
$u_i>0$. The other coordinates are nonnegative and independent.
Summing their geometric series, with the exponents in
\eqref{eq:orthogonalexponents}, gives \eqref{eq:orthogonalformula}.
The only non-root chamber facet is $u_0=0$. Its primitive form takes
value one on the central scalar cocharacter $(m;a)=(1;0,\ldots,0)$,
which has determinant degree $d$. The central translation condition,
split rank $r+1$, and longest length $r(r+1)$ give
\eqref{eq:orthogonalfe} by Theorem~\ref{thm:main}.
\end{proof}
Since the largest ray ratio is
$r(r+1)/(2r+2)=r/2$ and the numerator is positive at positive
$X,Y$, the abscissa of convergence is $r/2$. It is therefore
$1$, $3/2$, and $2$ in dimensions six, eight, and ten, respectively.

\subsection{Explicit formulas in dimensions six, eight, and ten}
In each formula below put $T=Y^d$, where $d$ is the indicated
dimension. Thus specialization means $X=q$ and $T=q^{-ds}$.

For $d=6$, the eight elements of $B_2$ give
\begin{equation}\label{eq:GSO6}
 \boxed{W_6^-(X,Y)=
 \frac{1+(1+2X+2X^3+X^4)T+X^4T^2}
      {(1-T)(1-X^4T)(1-X^6T)},\qquad T=Y^6.}
\end{equation}
For example its coefficient of $T$ at $X=q$ is
$q^6+2q^4+2q^3+2q+2$. Directly, the three possibilities
$a=(0,0),(1,0),(1,1)$ with $m=1$ have respective volumes
\[
 1,\qquad q^4+q^3+q+1,\qquad q^6+q^4+q^3+q,
\]
which give the same coefficient. The functional-equation factor
for \eqref{eq:GSO6} is $-X^6Y^6$.

For $d=8$, define the explicit polynomial
\begin{align*}
 A_8(X)={}&1+2X+2X^2+2X^3+2X^4+4X^5\\
          &\quad+3X^6+2X^7+2X^8+2X^9+X^{10}.
\end{align*}
The forty-eight elements of $B_3$ give
\begin{equation}\label{eq:GSO8}
 W_8^-(X,Y)=
 \frac{1+A_8(X)T+X^{16}A_8(X^{-1})T^2+X^{16}T^3}
      {(1-T)(1-X^6T)(1-X^{10}T)(1-X^{12}T)},\quad T=Y^8.
\end{equation}
The reflected expression $X^{16}A_8(X^{-1})$ is an ordinary
polynomial with exponents six through sixteen, not an additional
unknown coefficient. The functional-equation factor is $X^{12}Y^8$.

For $d=10$, let
\begin{align*}
 A_{10}(X)={}&1+2X+2X^2+4X^3+2X^4+4X^5+6X^6+6X^7+7X^8\\
 &+6X^9+6X^{10}+6X^{11}+6X^{12}+6X^{13}
      +3X^{14}+4X^{15}\\
 &+2X^{16}+2X^{17}+X^{18},\\
 C_{10}(X)={}&1+2X+4X^2+4X^3+6X^4+8X^5+9X^6+12X^7\\
 &+14X^8+14X^9+17X^{10}+16X^{11},\\
 B_{10}(X)={}&X^8\bigl(C_{10}(X)+16X^{12}
                            +X^{24}C_{10}(X^{-1})\bigr).
\end{align*}
The three hundred and eighty-four elements of $B_4$ give
\begin{equation}\label{eq:GSO10}
 \begin{split}
 W_{10}^-(X,Y)&=
 \frac{1+A_{10}(X)T+B_{10}(X)T^2
              +X^{40}A_{10}(X^{-1})T^3+X^{40}T^4}
      {(1-T)(1-X^8T)(1-X^{14}T)(1-X^{18}T)(1-X^{20}T)},\\
 &\hspace{7cm}T=Y^{10}.
 \end{split}
\end{equation}
Here too all numerator coefficients are explicit ordinary
polynomials. The functional-equation factor is $-X^{20}Y^{10}$.
All three numerator values at $X=T=1$ equal the respective Weyl
orders $8$, $48$, and $384$, an immediate additional check.

Although Tits states Example~1.16 for $r\geq2$, the same norm-plane
construction in rank one gives the useful dimension-four comparison
\[
 W_4^-(X,Y)=\frac{1+T}{(1-T)(1-X^2T)},\qquad T=Y^4,
\]
with functional-equation factor $X^2Y^4$. This is a lower-rank
extension of the calculation, rather than an additional case in
Tits's stated range.

\section{Ring realizations and pro-isomorphic subring zeta functions}
\label{sec:rings}
We now adapt the bilinear-form construction of
\cite[\S3.3]{duSL} to the quasi-split examples above. Throughout
this section $k=\Q_p$, $\OO=\Z_p$, and $p$ is odd. All rings
constructed below are finite free over $\OO$, associative and
nilpotent, with zero triple products; they need not have an identity
element and need not be Lie rings.

For clarity, the group pro-isomorphic zeta function is
\[
 \zeta^\wedge_\Gamma(s)=
 \sum_{\substack{H\leq\Gamma,\ [\Gamma:H]<\infty\\
                    \widehat H\cong\widehat\Gamma}}[\Gamma:H]^{-s}.
\]
Here we use its \emph{ring analogue}. For an $\OO$-ring $\mathcal A$
put
\begin{equation}\label{eq:ringzetadef}
 \zeta^\wedge_{\mathcal A}(s)=
 \sum_{\substack{M\leq_{\mathrm{ring}}\mathcal A,\ [\mathcal A:M]<\infty\\
                     M\cong_{\OO\text{-alg}}\mathcal A}}
                    [\mathcal A:M]^{-s}.
\end{equation}
These are isomorphisms of topological $\Z_p$-rings. For a finite
free $\Z$-ring the corresponding global condition is
$M\otimes\Z_p\cong\mathcal A\otimes\Z_p$ at every prime, and the
global series has the local factors \eqref{eq:ringzetadef}; see
\cite[Proposition~1.1]{duSL}. We do not identify these associative
ring zeta functions with zeta functions of nilpotent groups.

\subsection{A bilinear-form ring and its unipotent fibre}
Let $V=k^D$, let $V_{\OO}=\OO^D$, and let
$B:V\times V\to k$ be a bilinear form whose matrix is in
$\GL_D(\OO)$. Symmetry or alternation is not required. Define
\begin{equation}\label{eq:bilinearring}
 \mathcal A_B=V_{\OO}\oplus\OO z,\qquad
 (v+az)(w+bz)=B(v,w)z.
\end{equation}
All products involving $z$ vanish, so associativity follows from
$\mathcal A_B^3=0$. Its square and its two-sided annihilator are
both $\OO z$. This is the characteristic annihilator, rather than
the commutative centre, which in the symmetric case is the entire ring.
Write
\[
 H_B=\{g\in\GL(V):B(gv,gw)=c(g)B(v,w)
                          \text{ for some }c(g)\in k^\times\}.
\]
The forms used below have reductive identity component $H_B^\circ$.

\begin{proposition}\label{prop:ringautomorphisms}
Suppose that $H_B^\circ$ is reductive. For the $k$-algebra
$A_B=\mathcal A_B\otimes_{\OO}k$,
every automorphism has the form
\begin{equation}\label{eq:ringautmap}
 v+az\longmapsto gv+(\ell(v)+c(g)a)z,
 \qquad g\in H_B,\quad\ell\in\operatorname{Hom}_k(V,k).
\end{equation}
Thus the automorphism group is a split extension of $H_B$ by the
vector group $N=\operatorname{Hom}_k(V,k)$. Then $N$ is the
unipotent radical of the identity component,
and the full quotient is $H_B$.
With $\rho_{\mathrm{nat}}:H_B\hookrightarrow\GL_D$, one has
\begin{equation}\label{eq:ringshift}
 \boxed{\zeta^\wedge_{\mathcal A_B}(s)=
 Z_{H_B,\rho_{\mathrm{nat}}}
               \left(\left(1+\frac2D\right)s-2\right).}
\end{equation}
The measure on $H_B(k)$ is normalized at the integral automorphisms
of $V_{\OO}$.
\end{proposition}
\begin{proof}
The characteristic annihilator forces an automorphism to act on $z$
by a nonzero scalar $c$ and on $A_B/kz$ by some $g\in\GL(V)$.
Preservation of multiplication is exactly
$B(gv,gw)=cB(v,w)$; the linear functional $\ell$ is unrestricted.
Conversely every such map preserves multiplication. This proves
\eqref{eq:ringautmap}, the splitting, and the assertion about the
unipotent radical.

We give a lattice count to determine the fibre without a measure
convention ambiguity. For $g\in H_B$ with $gV_{\OO}\subseteq V_{\OO}$
put $U=gV_{\OO}$ and $c=c(g)$. Taking determinants in the
similitude identity gives
\begin{equation}\label{eq:multiplierdet}
 |c|=|\det g|^{2/D}.
\end{equation}
In particular $c\in\OO$. An isomorphic subring with projection $U$
has intersection $c\OO z$ with $kz$. It is specified by a graph
map $U\to\OO/c\OO$, and every such map occurs by choosing $\ell$
in \eqref{eq:ringautmap}. Since $U$ has rank $D$, the number of
graphs is $|c|^{-D}$. Distinct graphs give distinct subrings.
Their index is $|\det g|^{-1}|c|^{-1}$.

The possible projection lattices are indexed by the integral
similitudes modulo $H_B(\OO)$. The valuation of $c$ depends only
on the projection lattice by \eqref{eq:multiplierdet}; changing its
similitude basis changes $c$ by a unit. The normalized Haar integral
over the reductive group counts these lattices. Consequently
\begin{align}
 \zeta^\wedge_{\mathcal A_B}(s)
  &=\int_{H_B^+}|\det g|^s|c(g)|^{s-D}\,d\mu_{H_B}(g)\notag\\
  &=\int_{H_B^+}|\det g|^{(1+2/D)s-2}\,d\mu_{H_B}(g),
                                                        \label{eq:ringfibre}
\end{align}
which is \eqref{eq:ringshift}. This also verifies integral lifting:
the diagonal lift $(g,c(g))$ is integral whenever $g$ is integral.
\end{proof}

\begin{remark}[The exponent contributed by the fibre]
\label{rem:fibrenormalization}
The fibre factor is $|c|^{-D}=|\det g|^{-2}$, not
$|\det g|^{-D}$. Thus the constant in the change of variable is
$-2$. This corrects the factor written in
\cite[p.~64, Lemma~3.14]{duSL}: the $D$ additive parameters in
the fibre are scaled by the action on the one-dimensional
annihilator. An elementary check is the rank-two ring with
$x^2=z$ and all other basis products zero. Its isomorphic subrings are
\[
 \langle p^m x+uz,\ p^{2m}z\rangle_{\Z_p},\qquad
 m\geq0,\quad u\in\Z_p/p^{2m}\Z_p.
\]
Their index is $p^{3m}$ and their number for fixed $m$ is $p^{2m}$,
so the series is $(1-p^{2-3s})^{-1}$. This agrees with
\eqref{eq:ringshift} for $D=1$.
\end{remark}

\subsection{Realizing the quasi-split unitary similitude groups}
Let $E=k(\tau)$ be unramified quadratic, with
$\tau^2=\epsilon\in\OO^\times$ a nonsquare unit and
$\bar\tau=-\tau$. Let $J_n$ be the anti-diagonal $n\times n$
matrix with all anti-diagonal entries one, and put
\[
 h(x,y)=\bar x^{\,t}J_ny\quad(x,y\in E^n),\qquad
 \beta=1+\tau.
\]
This is the quasi-split hermitian form used in the preceding unitary
calculations. On the underlying $2n$-dimensional $k$-space define
\begin{equation}\label{eq:unitarybilinear}
 B_h(x,y)=\frac12\Tr_{E/k}(\beta h(x,y)),\qquad
 \mathcal U_n=\mathcal A_{B_h}.
\end{equation}
The ring $\mathcal U_n$ has $\Z_p$-rank $2n+1$.
Writing $x=u+\tau v$ and ordering the basis with the $u$-coordinates
first, an explicit integral multiplication matrix is
\begin{equation}\label{eq:unitaryringmatrix}
 B_h=\begin{pmatrix}J_n&\epsilon J_n\\
                  -\epsilon J_n&-\epsilon J_n\end{pmatrix}.
\end{equation}
Its determinant is $\epsilon^n(\epsilon-1)^n\det(J_n)^2$, a unit.

\begin{lemma}\label{lem:unitaryringquotient}
The $k$-linear similitude group of $B_h$ is exactly $\GU(h)$:
the condition of $E$-linearity follows from $B_h$ itself and need
not be imposed as extra structure. In particular
\[
 \operatorname{Aut}_k(\mathcal U_n\otimes k)/R_u
       \cong\GU(h),\qquad \dim R_u=2n.
\]
\end{lemma}
\begin{proof}
For a nonsingular bilinear form its asymmetry operator $C$ is
defined intrinsically by $B_h(x,Cy)=B_h(y,x)$. Here $C$ is
multiplication by $\bar\beta/\beta$. This element generates $E/k$:
it is neither $1$ nor $-1$, and an element of $k$ of norm one would
have to be one of these. Every similitude of $B_h$ commutes with
$C$, hence with the $k$-algebra $k[C]=E$; it is therefore $E$-linear.
For such a map $g$, preservation of $B_h$ with multiplier $c\in k^\times$
means
\[
 \Tr\bigl(\beta(h(gx,gy)-c h(x,y))\bigr)=0.
\]
Replacing $y$ by $\tau y$ gives the same equation with $\beta\tau$
in place of $\beta$. Nondegeneracy of the trace pairing forces
$h(gx,gy)=c h(x,y)$. Conversely an $E$-linear unitary similitude
preserves \eqref{eq:unitarybilinear} with the same multiplier.
The assertion follows from Proposition~\ref{prop:ringautomorphisms}.
\end{proof}
The nonsymmetric choice of $B_h$ is essential to this realization:
its asymmetry records the quadratic extension and excludes unwanted
semilinear transformations. A symmetric trace form alone would give
an orthogonal similitude group on the underlying $k$-space.

\begin{theorem}\label{thm:unitaryringzeta}
For the ring $\mathcal U_n$,
\begin{equation}\label{eq:unitaryringshift}
 \zeta^\wedge_{\mathcal U_n}(s)=
 Z_{\GU(n),\rho_{\mathrm{nat}}}
               \left(\frac{n+1}{n}s-2\right),
\end{equation}
where $\rho_{\mathrm{nat}}$ is the $2n$-dimensional $k$-representation
used in Sections~\ref{sec:GU3}--\ref{sec:GUn}.
Thus all the preceding finite Weyl-group formulas give pro-isomorphic
subring zeta functions by this substitution.
\end{theorem}
\begin{proof}
Apply Proposition~\ref{prop:ringautomorphisms} with $D=2n$ and
Lemma~\ref{lem:unitaryringquotient}. The action on the entire ring
has determinant $\det_k(g)c(g)=c(g)^{n+1}$, whereas the unipotent
fibre has measure $|c(g)|^{-2n}$. These give precisely the exponent
in \eqref{eq:unitaryringshift}.
\end{proof}

For example, in dimension $n=2$ the ring has rank five and
\[
 \zeta^\wedge_{\mathcal U_2}(s)=
 \frac1{(1-p^{4-3s})(1-p^{5-3s})}.
\]
For $n=3$, the rank-seven ring has
\begin{equation}\label{eq:unitaryring3}
 \boxed{\zeta^\wedge_{\mathcal U_3}(s)=
 \frac{1+p^{13-8s}}{(1-p^{12-8s})(1-p^{16-8s})}.}
\end{equation}
For $n=4$, the rank-nine ring has
\begin{equation}\label{eq:unitaryring4}
 \zeta^\wedge_{\mathcal U_4}(s)=
 \frac{(1+p^{9-5s})(1+p^{11-5s})}
      {(1-p^{8-5s})(1-p^{12-5s})(1-p^{22-10s})}.
\end{equation}
For $n=5$, put $T=p^{20-12s}$. The rank-eleven ring has
\[
 \zeta^\wedge_{\mathcal U_5}(s)=
 \frac{1+(p+p^3+p^4+p^6+p^7+p^9)T+p^{10}T^2}
      {(1-T)(1-p^8T)(1-p^{12}T)}.
\]

If the quadratic algebra is split, the same bilinear construction
works with the split involution, provided $\beta$ and the difference
of the two asymmetry eigenvalues are units. It realizes the split
unitary similitude group. In particular, for the rank-seven ring at
a good split prime, putting $U=p^{6-4s}$ gives
\begin{equation}\label{eq:unitaryring3split}
 \zeta^\wedge_{\mathcal U_3,\mathrm{sp}}(s)=
 \frac{1+2(p+1)U+pU^2}{(1-U)^2(1-p^2U)^2}.
\end{equation}
This differs from \eqref{eq:unitaryring3} at nonsplit primes.

\subsection{Realizing Tits's nonsplit orthogonal forms}
Let $Q_r$ be \eqref{eq:orthogonalform}, on a space of dimension
$d=2r+2$. Take its polar bilinear form
\[
 B_Q(v,w)=Q_r(v+w)-Q_r(v)-Q_r(w)
\]
and define $\mathcal O_d=\mathcal A_{B_Q}$. This is a commutative
associative ring of rank $d+1$. In paired hyperbolic coordinates and
the norm-plane basis $1,\tau$, its multiplication matrix has blocks
\[
 \begin{pmatrix}0&1\\1&0\end{pmatrix}\quad(r\text{ times}),
 \qquad\begin{pmatrix}2&0\\0&-2\epsilon\end{pmatrix}.
\]
This matrix is unimodular at the primes under consideration.
Since $p$ is odd, $H_{B_Q}=\GO(Q_r)$. Hence
\begin{align*}
 \operatorname{Aut}_k(\mathcal O_d\otimes k)/R_u
      &\cong\GO(Q_r),\qquad \dim R_u=d,\\
 \operatorname{Aut}_k(\mathcal O_d\otimes k)^\circ/R_u
      &\cong\GSO(Q_r).
\end{align*}
The full quotient has two components; it is the connected quotient
that is $\GSO$. Conjugation on the norm plane is an integral
improper isometry, so the compact subgroup meets both components.
As in Section~\ref{sec:orthogonal}, full and connected normalized
integrals are equal, with no additional factor of two.

\begin{theorem}\label{thm:orthogonalringzeta}
For $d=2r+2$ and the unramified nonsplit form of Tits's
Example~1.16,
\begin{equation}\label{eq:orthogonalringshift}
 \zeta^\wedge_{\mathcal O_d}(s)=
 Z_{\GSO(Q_r),\rho_{\mathrm{nat}}}
                 \left(\frac{d+2}{d}s-2\right).
\end{equation}
More explicitly, let $L$ on $B_r$ have parameters $(1,\ldots,1,2)$,
put $T=p^{2d-(d+2)s}$, and set $M_i=p^{i(2r+1-i)}T$. Then
\begin{equation}\label{eq:orthogonalringweyl}
 \zeta^\wedge_{\mathcal O_d}(s)=
 \frac{\displaystyle\sum_{w\in B_r}p^{-L(w)}
                                  \prod_{i\in\Des(w)}M_i}
      {(1-T)\displaystyle\prod_{i=1}^r(1-M_i)}.
\end{equation}
\end{theorem}
\begin{proof}
Use Proposition~\ref{prop:ringautomorphisms} and the compact
representative of the improper component to replace $\GO$ by $\GSO$.
Substitution in \eqref{eq:orthogonalformula} gives
\eqref{eq:orthogonalringweyl}.
\end{proof}

For $d=6$, the resulting rank-seven ring has, with $T=p^{12-8s}$,
\begin{equation}\label{eq:orthogonalring6}
 \boxed{\zeta^\wedge_{\mathcal O_6}(s)=
 \frac{1+(1+2p+2p^3+p^4)T+p^4T^2}
      {(1-T)(1-p^4T)(1-p^6T)}.}
\end{equation}
For $d=8$, the rank-nine ring has, with $T=p^{16-10s}$,
\begin{equation}\label{eq:orthogonalring8}
 \zeta^\wedge_{\mathcal O_8}(s)=
 \frac{1+A_8(p)T+p^{16}A_8(p^{-1})T^2+p^{16}T^3}
      {(1-T)(1-p^6T)(1-p^{10}T)(1-p^{12}T)}.
\end{equation}
For $d=10$, the rank-eleven ring has, with $T=p^{20-12s}$,
\begin{equation}\label{eq:orthogonalring10}
 \zeta^\wedge_{\mathcal O_{10}}(s)=
 \frac{1+A_{10}(p)T+B_{10}(p)T^2+
                     p^{40}A_{10}(p^{-1})T^3+p^{40}T^4}
      {(1-T)(1-p^8T)(1-p^{14}T)(1-p^{18}T)(1-p^{20}T)}.
\end{equation}
The polynomials $A_8,A_{10},B_{10}$ are explicitly listed in
Section~\ref{sec:orthogonal}; thus these are fully specified rational
functions. The abscissae in dimensions $6,8,10$ are respectively
$9/4$, $14/5$, and $10/3$.

\subsection{Functional equations, global rings, and scope}
The changes of variable above also transport the functional equations.
Simultaneous inversion means $p\mapsto p^{-1}$ and $p^{-s}\mapsto p^s$.
For the unitary rings at nonsplit primes, with $r=\lfloor n/2\rfloor$,
\begin{equation}\label{eq:unitaryringfe}
 \left.\zeta^\wedge_{\mathcal U_n}(s)\right|_{p\to p^{-1}}
  =(-1)^{r+1}p^{n(n-1)/2+4n-2(n+1)s}
                      \zeta^\wedge_{\mathcal U_n}(s).
\end{equation}
At split primes the same formula holds with sign $(-1)^{n+1}$.
For the orthogonal rings, $d=2r+2$,
\begin{equation}\label{eq:orthogonalringfe}
 \left.\zeta^\wedge_{\mathcal O_d}(s)\right|_{p\to p^{-1}}
  =(-1)^{r+1}p^{r(r+1)+2d-(d+2)s}
                      \zeta^\wedge_{\mathcal O_d}(s).
\end{equation}
For instance, the rank-seven unitary and orthogonal rings have
factors $p^{15-8s}$ and $-p^{18-8s}$ respectively.

These local examples can be obtained from explicit $\Z$-rings.
Fix a nonsquare integer $\epsilon\ne0,1$ and use the integer matrix
\eqref{eq:unitaryringmatrix} for the unitary ring and the displayed
polar matrix for the orthogonal ring. Away from the primes dividing
$2\epsilon(\epsilon-1)$, the unitary reduction is governed by the
quadratic \emph{\'etale} algebra $\Q_p[\tau]/(\tau^2-\epsilon)$.
At nonsplit primes the extension is unramified and the formulas above
apply; at split primes use the split group formulas, with the same
change of variable. In particular \eqref{eq:unitaryring3} and
\eqref{eq:unitaryring3split} give the two good local factors of one
rank-seven integer ring. For the orthogonal integer ring the same
construction gives the displayed nonsplit factors at good inert
primes; its split factors are the corresponding split orthogonal
similitude integrals under \eqref{eq:ringshift}. Ramified and other
excluded primes are not covered by these good-prime formulas.

The constructions are ring realizations in the precise sense of
\cite[\S3.3]{duSL}. They do not themselves produce nilpotent groups
with these pro-isomorphic factors via the Malcev correspondence:
the orthogonal multiplication is symmetric, and the unitary
multiplication is not alternating. Taking commutators would discard
part of the multiplication. In the orthogonal case it gives the zero
Lie bracket; in the unitary case it retains only a nondegenerate
alternating form and hence a symplectic, rather than unitary,
similitude group. This explains why the ring and group interpretations
must be kept distinct here.

\end{document}